\documentclass{sig-alternate}

\begin{document}

\permission{Author's version.}

\conferenceinfo{GECCO'15,} {July 11-15, 2015, Madrid, Spain.}

    \clubpenalty=10000
    \widowpenalty = 10000

\title{A Closer Look At Differential Evolution For The Optimal Well Placement Problem}
    
\numberofauthors{4}

\author{
\alignauthor
Grazieli L. C. Carosio\titlenote{Corresponding author.}\\
       \affaddr{Department of Mathematics and Statistics}\\
       \affaddr{Memorial University of Newfoundland}\\
       \affaddr{St. John's, Canada A1C 5S7}\\
       \email{grazielilccdm@mun.ca}
\alignauthor
Thomas D. Humphries\\
       \affaddr{Department of Mathematics}\\
       \affaddr{Oregon State University}\\
       \affaddr{Corvallis, USA 97331}\\
       \email{humphrit@science.oregonstate.edu}
\and
\alignauthor Ronald D. Haynes\\
       \affaddr{Department of Mathematics and Statistics}\\
       \affaddr{Memorial University of Newfoundland}\\
       \affaddr{St. John's, Canada A1C 5S7}\\
       \email{rhaynes@mun.ca}
\\  
\alignauthor Colin G. Farquharson\\
       \affaddr{Department of Earth Sciences}\\
       \affaddr{Memorial University of Newfoundland}\\
       \affaddr{St. John's, Canada A1B 3X5}\\
       \email{cgfarquh@mun.ca}
       }

\maketitle
\begin{abstract}
Energy demand has increased considerably with the growth of world population, increasing the interest in the hydrocarbon reservoir management problem. Companies are concerned with maximizing oil recovery while minimizing capital investment and operational costs. A first step in solving this problem is to consider optimal well placement. In this work, we investigate the Differential Evolution (DE) optimization method, using distinct configurations with respect to population size, mutation factor, crossover probability, and mutation strategy, to solve the well placement problem. By assuming a bare control procedure, one optimizes the parameters representing positions of injection and production wells. The Tenth SPE Comparative Solution Project and MATLAB Reservoir Simulation Toolbox (MRST) are the benchmark dataset and simulator used, respectively. The goal is to evaluate the performance of DE in solving this important real-world problem. We show that DE can find high-quality solutions, when compared with a reference from the literature, and a preliminary analysis on the results of multiple experiments gives useful information on how DE configuration impacts its performance.
\end{abstract}




\section{Introduction}

Recently, energy demand has increased considerably with the growth of world population, inducing even more interest in the oil reservoir management problem. This involves effectively utilizing existing hydrocarbon resources by maximizing recovery and minimizing capital investment and operational costs. One approach to address this is to find optimal placement of wells. The parameters representing positions and orientations of injection and production wells are optimized, while a bare control procedure is assumed. In this work, we consider only vertical wells with $(x,y)$-coordinates representing a well's position.

The performance indicator or objective function of this optimization problem is usually either the Voidage Replacement Ratio (VRR)~\cite{Clark2003,Awotunde2014a} or the Net Present Value (NPV)~\cite{Bangerth2006,Humphries2014}. The VRR function corresponds to the ratio of the total volume of fluid injected (water injection) to the volume of fluid produced (water, oil and gas production). In turn, the NPV function involves the total amount of oil extracted, emphasizing an early oil production in the reservoir's life-time (due to the time value of money), and also typically includes the costs of water injection and disposal of any water produced. Since NPV is the primary objective function used in well placement optimization, we chose it for this study. 

Since the reservoir model is represented as a discretized grid, the positions of injection and production wells are typically specified by integer co-ordinates, meaning that exact gradients of the objective function with respect to these parameters are not available. Thus, gradient-based methods have been used only after conversion of the discrete optimization problem into another optimization problem with continuous variables~\cite{Sarma2008,Zhang2010}. Due mostly to the heterogeneous permeability field, the objective function is non-convex and may contain many local minima~\cite{Onwunalu2010}. Moreover, gradient-based methods are typically local solvers. Here a global optimum is desired, hence the work presented herein uses a derivative-free (non-invasive, black-box) optimization method known as a meta-heuristic~\cite{Talbi2009,Boussaid2013}. 

Since gradients are not directly available, meta-heuristics are usually quite inefficient requiring several thousands of simulations and thus may have limited application to large-scale simulation models with many wells. However, whenever they are computationally feasible, meta-heuristics tend to provide high-quality solutions. Parallel programming can be used to obtain multiple independent runs simultaneously, or provide parallel function evaluations in a single run, thus the computational time may be decreased and the overall solution can be improved. Another important advantage of applying meta-heuristics is that the number of function evaluations is dependent on the size of the population and not the problem dimensionality~\cite{Temizel2014}. Meta-heuristics that have been used to solve the well-placement problem include: genetic algorithms (GA)~\cite{Yeten2003,Morales2010,Lyons2013}, simulated annealing (SA)~\cite{Bangerth2006}, particle swarm optimization (PSO)~\cite{Onwunalu2010,Isebor2013,Afshari2014,Humphries2014}, covariance matrix adaptation evolutionary strategies (CMA-ES)~\cite{Ding2008,Bouzarkouna2012,Fonseca2013}, and differential evolution (DE)~\cite{Nwankwor2013,Awotunde2014a}.

DE is a meta-heuristic which has been successfully used to solve many global optimization problems~\cite{Chakraborty2008}, for example, problems in engineering design~\cite{Melo2013a}, electrical impedance tomography~\cite{Leskinen2009}, the paper industry~\cite{Tirronen2008}, and limited memory optimization~\cite{Neri2011}. Nevertheless, its application in well-placement is not common. Some researchers~\cite{Nwankwor2013} indicate that DE is not appropriate for the well placement problem, suggesting its use only when hybridized with other stochastic methods. However, it is well-known that the choice of algorithmic parameters may have great influence on its behaviour. Therefore, in this paper we investigate DE to solve a well-placement problem on a widely studied real-world example, and present a preliminary analysis on its parameters. Our results may serve as reference for researchers intending to use DE to solve this kind of problem; allowing a more fair comparison by using a more appropriate DE configuration. Moreover, this study suggests a particular configuration that may make classical DE competitive with other methods to solve the well placement problem.

This paper is organized as follows: Section 2 presents the literature review; Section 3 presents the proposed methodology; Section 4 presents the results and discussions; Section 5 presents the conclusions.
~\\
\section{Literature Review}

The well placement problem involves determining the optimal location of one or more wells in order to obtain, for example, the maximum NPV of oil extracted with minimal operating costs, subject to different geological and economic constraints such as reservoir geology, well type, and physical properties of the fluid and rock. We focus this short review on DE, the optimization method investigated in this paper.

DE has been used in multi-objective approaches incorporating NPV and VRR to solve well placement optimization problems~\cite{Ahmed2014,Awotunde2014a}. A Pareto-based multi-objective approach to integrate NPV and VRR in a well placement optimization was presented by~\cite{Ahmed2014}, obtaining several optimal solutions according to the Pareto front. The set of solutions in the Pareto optimal front are those that cannot be improved in any of the objectives without degrading at least one of the other objectives.

A similar study was presented by~\cite{Awotunde2014a} using CMA-ES and DE to solve the multi-objective optimization problem by manually weighting NPV and VRR, resulting in a single-optimization problem. They concluded that CMA-ES outperformed DE for the majority of experiments, but that DE also obtained high-quality solutions.

A hybrid approach (HPSDE) was proposed by~\cite{Nwankwor2013} combining the DE and PSO algorithms. The hybrid algorithm consistently outperformed both DE and PSO using NPV as performance measure in the three proposed example problems. The placement of one, two, and nine vertical wells were considered. They showed that the performance of DE and PSO, to an extent, depends on the 1) the total number of simulations, and 2) the population size, and that PSO was generally better than DE. However, the conclusion of that study was based on results obtained from a configuration that makes DE's convergence slow but able to avoid local optima. Therefore, more iterations would be necessary to show a competitive performance with that configuration. Their conclusions were based only $5$ trials of each algorithm. We understand that an adequate statistical analysis of a stochastic algorithm requires a larger number of independent runs.

Given the above, our main contribution is to present a preliminary analysis using various DE configurations, including some suggested in the literature~\cite{Storn1997,Das2011}, while solving a well placement problem. We show that DE can find high-quality solutions, contradicting the conclusions pres\-ented by \cite{Nwankwor2013}, and indeed our DE obtained several results better than those from~\cite{Humphries2014}, which used a PSO with local-search (memetic PSO).
~\\
\section{Methodology}

In this section we present the well placement problem used as the case study, and introduce the Differential Evolution algorithm used to solve it.
~\\
\subsection{Case Study\label{case_study}}

\begin{figure}[h]
\centering
\includegraphics[width=0.7\columnwidth]{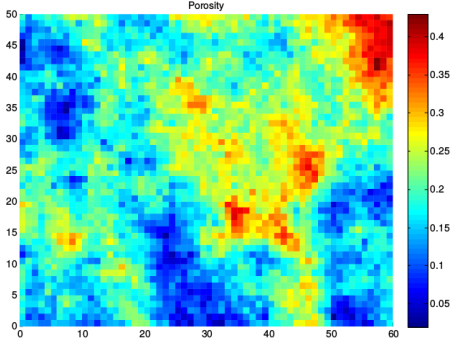} 
\caption{\label{poro}Porosity field.}
\end{figure}

\begin{figure}[h]
\centering
\includegraphics[width=0.7\columnwidth]{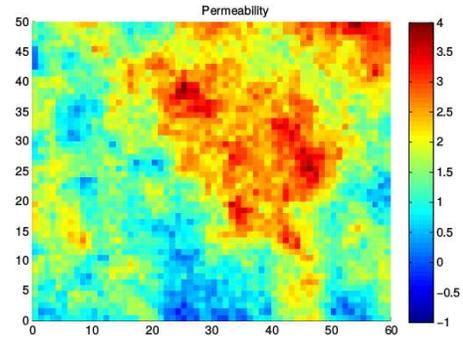} 
\caption{\label{perm}Permeability (mD) in logarithmic color scale.}
\end{figure}

As in~\cite{Humphries2014}, we consider a small section of the reservoir block representing a subset of the Tenth SPE Comparative Solution Project~\cite{Christie2001} to run several instances of each DE configuration and to study the convergence behaviour. We adopt a $2D$ reservoir model consisting of $60\times50$ grid cells, each measuring $32\times32$ m$^{2}$ with a total field size of $1,600\times1,920$ m$^{2}$. A uniform saturation of 80:20, oil to water, is also assumed as in~\cite{Humphries2014}. The physical parameters of the problem are shown in Table~\ref{parameters}. Figures~\ref{poro} and~\ref{perm} show the porosity and permeability fields of the reservoir model, which are generated from the third layer of the SPE10 base case.

\begin{table}[ht]
\begin{centering}
\caption{\label{parameters}Physical parameters values as used in~\cite{Humphries2014}.}
\par\end{centering}
\centering{}%
\begin{tabular}{|c|c|} \hline 
\textbf{Parameter} & \textbf{Value} \\ \hline 
Fluid viscosities $\mu_{o}$ and $\mu_{w}$  & $2.4$ and $1.0$ cp \\ \hline 
Fluid densities $\rho_{o}$ and $\rho_{w}$ & $835$ and $1,000$ kg/$m^{3}$ \\ \hline 
Initial reservoir pressure & $260$ bars \\ \hline 
Injector BHP range & $275-450$ bars \\ \hline 
Producer BHP range & $100-250$ bars \\ \hline 
Production period & $10$ years \\ \hline 
Control interval & $2$ years \\ \hline 
\end{tabular}
\end{table}

The objective function adopted is the NPV over the entire production period $[0,T]$. The NPV is computed as in~\cite{Bangerth2006}:
\begin{equation}
\begin{split}
NPV(\textrm{\textbf{x}}) = \int_{0}^{T}{\sum_{n \in prod}c_{o}q^{-}_{n,o}(\textrm{\textbf{x}},t)-c_{w,disp}q^{-}_{n,w}(\textrm{\textbf{x}},t)}\\ -\sum_{n \in inj}c_{w,inj}q^{+}_{n,w}(\textrm{\textbf{x}},t)(1+r)^{t}dt,
\end{split}
\end{equation}

\noindent where $\textrm{\textbf{x}}$ is the solution vector encoding the spatial coordinates of each well, for instance, $\textrm{\textbf{x}} = [x_{1},x_{2},...,x_{n}]$, where $(x_{1},x_{2})$ is the coordinate pair of the first well, $(x_{3},x_{4})$ represents the coordinate pair of the second well, and so on. The parameters $c_{o}$, $c_{w,disp}$, and $c_{w,inj}$ represent, respectively, the value of each barrel of oil produced, the cost of disposing of each barrel of produced water, and the cost of injecting a barrel of water into the reservoir. The functions $q_{n,o}^{-}(\textrm{\textbf{x}},t)$ and $q_{n,w}^{-}(\textrm{\textbf{x}},t)$ are the oil and water production rates (barrels/day) at a producer, while $q_{n,w}^{+}(\textrm{\textbf{x}},t)$ is the daily water injection rate at an injector. Note that these rates are determined from a reservoir simulation according to the value of $\textrm{\textbf{x}}$. The yearly interest rate is specified by $r$. The economic parameter values shown in Table~\ref{economic} are the same as provided in~\cite{Humphries2014} for all experiments.

\begin{table}[ht]
\begin{centering}
\caption{\label{economic}Economic parameters values as in~\cite{Humphries2014}.}
\par\end{centering}
\centering{}%
\begin{tabular}{|c|c|} \hline 
\textbf{Parameter} & \textbf{Value} \\ \hline 
$c_{o}$ & $80/$bbl \\ \hline 
$c_{w,disp}$ & $12/$bbl \\ \hline 
$c_{w,inj}$ & $8/$bbl \\ \hline 
$r$ & $10\%$ or $0\%$ \\ \hline 
Threshold for water cut & $78\%$ \\ \hline 
\end{tabular}
\end{table}

The experiments to calculate the NPV were completed using the Matlab Reservoir Simulation Toolbox (MRST)~\cite{Lie2012} which solves the flow and transport equations in alternating steps, in order to determine the phase pressures, flow rates and saturation at every time point. Modelling of simple vertical wells is provided using the Peaceman model~\cite{Peaceman1978}.

The optimization problem consisted of the placement of two injectors and two producers in a reservoir model under conditions according to Table~\ref{parameters}, constrained to a minimum distance between wells equal to $250$m.

We considered three optimization subproblems as in~\cite{Humphries2014} in order to compare those results with the ones obtained in this paper:
\begin{itemize}
\item Case 1: no constraints regarding injection and production rates, and $r=10\%$;
\item Case 2: no constraints regarding injection and production rates, and $r=0\%$;
\item Case 3: maximum flow constraints regarding to injection and production rates equal to $1,000$ m$^{3}$/day, and $r=10\%$.  
\end{itemize} 

The next section presents the meta-heuristic employed here to solve the well-placement problem.
~\\
\subsection{Differential Evolution\label{DE}}

DE~\cite{Storn1997} is a floating-point encoding populational meta-heuristic, working as an evolutionary algorithm, with only a few control parameters \cite{Das2011,Boussaid2013}. Several papers have shown that the classical DE and its variants outperform many other optimization methods, in terms of both convergence speed and robustness, when applied to hard benchmark functions and real-world problems. These problems include unconstrained global optimization~\cite{Zou2013}, constrained optimization \cite{Melo2013a}, multi-objec\-tive optimization \cite{Mezura-Montes2008}, large-scale optimization~\cite{Brest2010,Zhou2012}, and optimization in dynamic and uncertain environments~\cite{Mendes2005,duPlessis2012}. 

The DE algorithm proceeds by randomly initializing (commonly using a uniform distribution) a population of $D$-di\-mensional vectors inside the problem bounds, and evaluating the objective/fitness function for all the vectors in the population. Then, until a stopping condition is satisfied, the algorithm performs an iterative evolutionary process of mutation, crossover, and selection operators.

For each vector $\textrm{\textbf{x}}_{i}$ in a population of size $N$, the mutation operator uses the weighted difference of parent solutions to generate mutation vectors $\textrm{\textbf{v}}_{i}$. The two well-known mutation strategies investigated in this work are \textbf{rand/1} and \textbf{current-to-best/1}~\cite{Neri2010}, which are represented, respectively, by Eqs.~\ref{rand} and \ref{current}:
\begin{equation}
\label{rand}
\textrm{\textbf{v}}_{i}=\textrm{\textbf{x}}_{r1}+F.(\textrm{\textbf{x}}_{r2}-\textrm{\textbf{x}}_{r3}),
\end{equation}
\begin{equation}
\label{current}
\textrm{\textbf{v}}_{i}=\textrm{\textbf{x}}_{i}+F.(\textrm{\textbf{x}}_{best}-\textrm{\textbf{x}}_{i})+F.(\textrm{\textbf{x}}_{r1}-\textrm{\textbf{x}}_{r2}),
\end{equation}

\noindent where $\textrm{\textbf{x}}_{r1}$, $\textrm{\textbf{x}}_{r2}$ and $\textrm{\textbf{x}}_{r3}$ are three distinct and randomly chosen vectors from the population, $\textrm{\textbf{x}}_{best}$ is the best vector from the population, $F\in[0,2]$ is the mutation factor, and $.$ denotes a scalar-vector product. In the \textit{rand/1} strategy, the base vector to be perturbed is randomly selected from the population, and will move in the direction of the other two random vectors. Therefore, it is possible that only poor solutions are chosen to generate offspring. Obviously, the repetition of this could result in a poor search. Given that the best solution of the population is eventually chosen, the search tends to the best solution through the optimization process~\cite{Islam2012}. However, it is also clear that the bigger the population the lower the chance of selecting the current best solution.

On the other hand, the \textit{current-to-best/1} strategy uses as a base vector each member of the population; thus all best solutions will always be chosen. Moreover, the best solution found so far is always used to calculate the difference vector, not only guiding the search towards the current best but also performing a localized search. In this case, smaller populations give faster convergence with lower exploratory capability and may result in the search getting trapped in low-quality optima~\cite{Islam2012}.

The crossover operator is applied on $\textrm{\textbf{v}}_{i}$ to generate the final offspring vector $\textrm{\textbf{u}}_{i}$ whose $j$-th component is given as
\begin{equation}
\label{crossover}
\textrm{\textbf{u}}_{i,j}= \left\{ \begin{array}{ll} \textrm{\textbf{v}}_{i,j} & \mbox{if $U(0,1) \leq CR$ or $j=j_{rand}$},\\
\textrm{\textbf{x}}_{i,j} & \mbox{otherwise},\end{array} \right.
\end{equation}

\noindent where $U(0,1)$ is a random floating-point number from a uniform distribution between $0$ and $1$ generated for each $j$, $CR\in[0,\,1]$ is the crossover probability, and $j_{rand}$ is a randomly chosen index from $[1,D]$.

Finally, the selection operator selects the best evaluated vector between $\textrm{\textbf{x}}_{i}$ and $\textrm{\textbf{u}}_{i}$. The offspring replaces the parent if its fitness value is better. Otherwise, the parent is maintained in the population.

As one may notice, DE has four parameters: $N$, $F$, $CR$, and the mutation strategy. Therefore, in this paper we experiment with several configurations to show that, when properly configured, DE can have a good performance in a well placement problem.
~\\
\section{Experimental Analysis}

In this paper we investigate the performance of DE using eight distinct configurations to solve the three well-placement problems described in Section~\ref{case_study}. The distinct configurations are specified in Table~\ref{confs}, and suggested by~\cite{Zhang2009,Das2011}. The stopping criterion is 10,000 objective function evaluations (simulations using MRST) as in~\cite{Humphries2014}.

\begin{table}[ht]
\begin{centering}
\caption{\label{confs}DE configurations.}
\par\end{centering}
\centering{}
\begin{tabular}{|c|c|c|c|c|}
\hline 
\textbf{Configuration} & \textbf{$N$} & \textbf{$CR$} & \textbf{$F$} & \textbf{Mutation}\tabularnewline
\hline 
$1$ & $100$ & $0.5$ & $0.9$ & rand/1 \tabularnewline
\hline 
$2$ & $100$ & $0.9$ & $0.5$ & rand/1 \tabularnewline
\hline 
$3$ & $200$ & $0.5$ & $0.9$ & rand/1 \tabularnewline
\hline 
$4$ & $200$ & $0.9$ & $0.5$ & rand/1 \tabularnewline
\hline 
$5$ & $100$ & $0.5$ & $0.9$ & current-to-best/1 \tabularnewline
\hline
$6$ & $100$ & $0.9$ & $0.5$ & current-to-best/1 \tabularnewline
\hline 
$7$ & $200$ & $0.5$ & $0.9$ & current-to-best/1 \tabularnewline
\hline 
$8$ & $200$ & $0.9$ & $0.5$ & current-to-best/1 \tabularnewline
\hline  
\end{tabular}
\end{table}

As can be observed in Table~\ref{confs}, four parameters were modified. More extensive tests, using a larger number of configurations, could not be performed because the simulation is computationally expensive. Therefore, it is possible that an untested configuration leads to even better results than those in the current work. However, those eight configurations allow the identification of interesting behaviours.

Since meta-heuristics have stochastic components, it is necessary to perform multiple runs in order to assess the average performance. Thus, each configuration was independently run $30$ times, with distinct seeds, to allow for a more accurate comparison.

A death penalty approach was adopted for constraint handling. That means 1) if DE generates a solution for well placement in which the minimum distance between wells is less than $250$m, then the NPV for this solution is set to minus infinity; and 2) if the maximum flow of a well during the production time is above $1,000$ m$^{3}$/day, then the simulation is interrupted and the NPV is also set to minus infinity (only applied in Case~3).
~\\
\subsection{Results}

Tables~\ref{results1}-\ref{results3} show the results for each DE configuration for each case described in Section~\ref{case_study} with respect to fitness value, i.e, NPV. The first column of the tables labels the methods (DE configurations and memetic PSO~\cite{Humphries2014}). The second column gives the best NPV found during the optimization process, while the third shows the worst NPV obtained. The mean, standard deviation (SD), and median values of NPV are shown in the fourth, fifth, and sixth columns respectively. The last line of the tables presents the results from~\cite{Humphries2014}, which were obtained by decoupling the well placement and control problems. That approach can be described in three steps: 1) PSO was used to determine the optimal well positions under a fixed control procedure; 2) once optimal positions were found under the fixed control procedure, the controls and positions were optimized locally using general pattern search (GPS) with standard search directions. To generate our results with DE we used the optimal control strategy found in~\cite{Humphries2014} and then determined optimal well placement.

\begin{figure}[ht]
\centering
\includegraphics[width=\columnwidth]{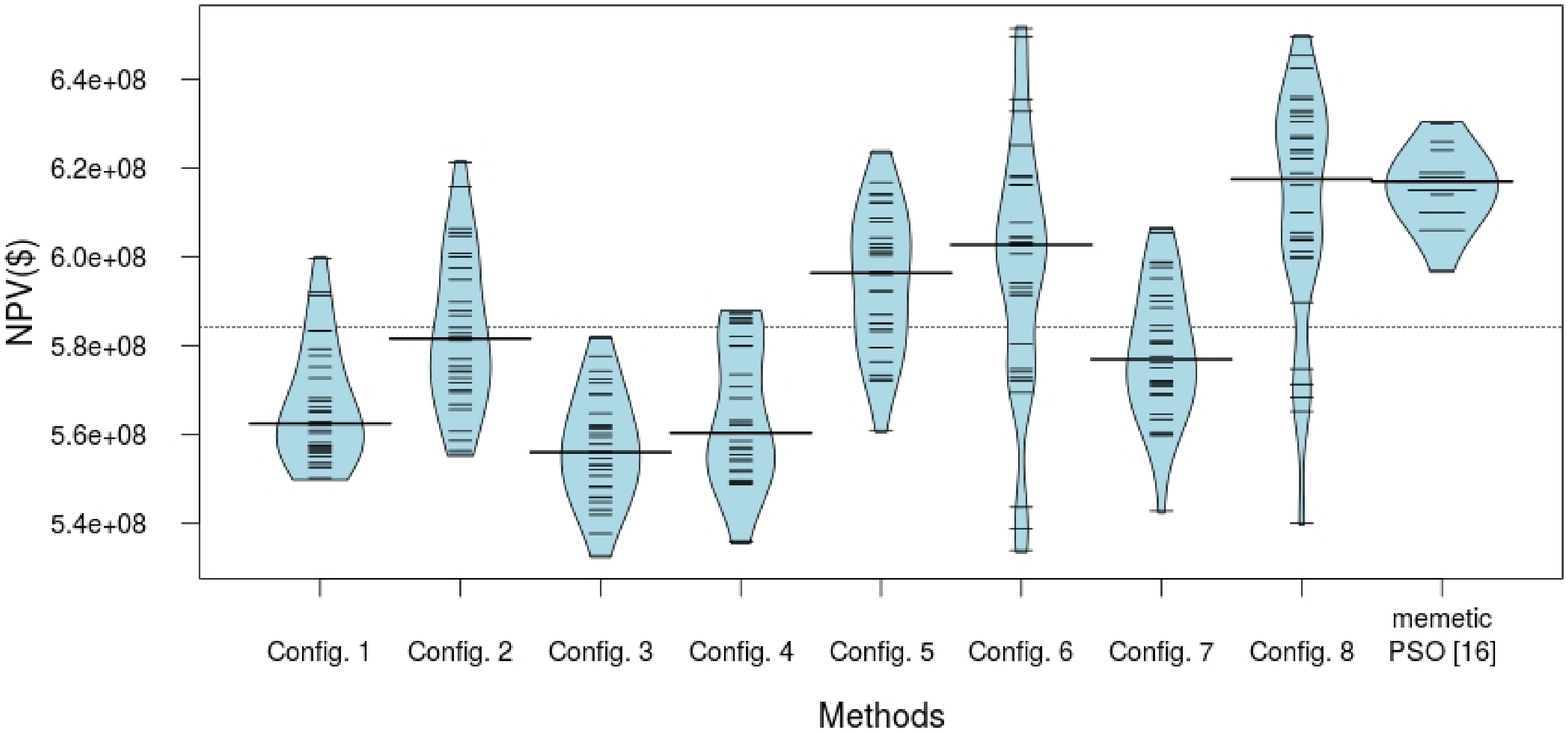} 
\caption{\label{bean1}Bean plots of the NPV distribution for Case~1. The short horizontal lines are the NPVs for all trials. The median is represented by the thick horizontal black line. The average of all observations is presented as the dashed horizontal line.}
\end{figure}

\begin{figure}[h]
\centering
\includegraphics[width=\columnwidth]{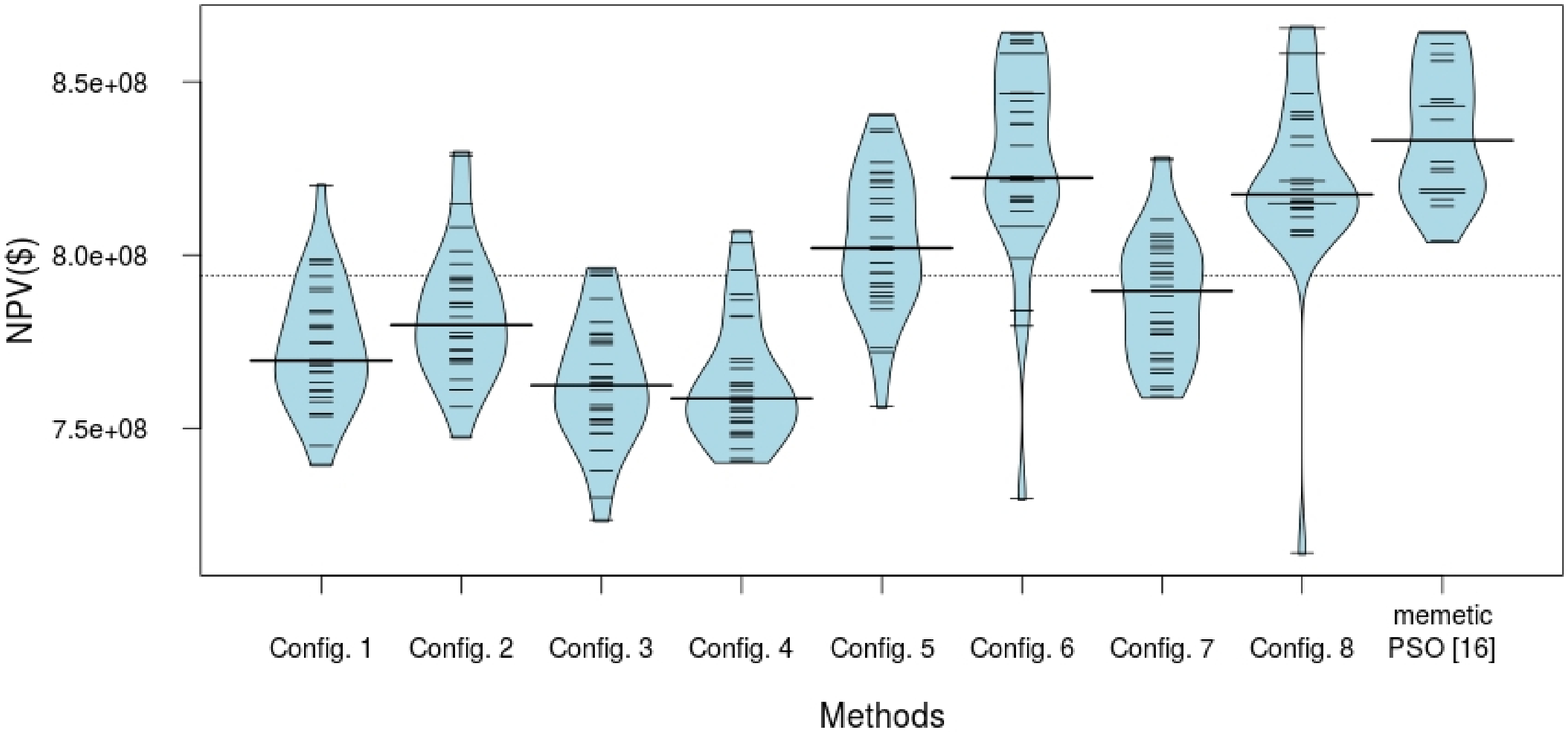} 
\caption{\label{bean2}Bean plots of the NPV distribution for Case~2. The short horizontal lines are the NPVs for all trials. The median is represented by the thick horizontal black line. The average of all observations is presented as the dashed horizontal line.}
\end{figure}

\begin{figure}[h]
\centering
\includegraphics[width=\columnwidth]{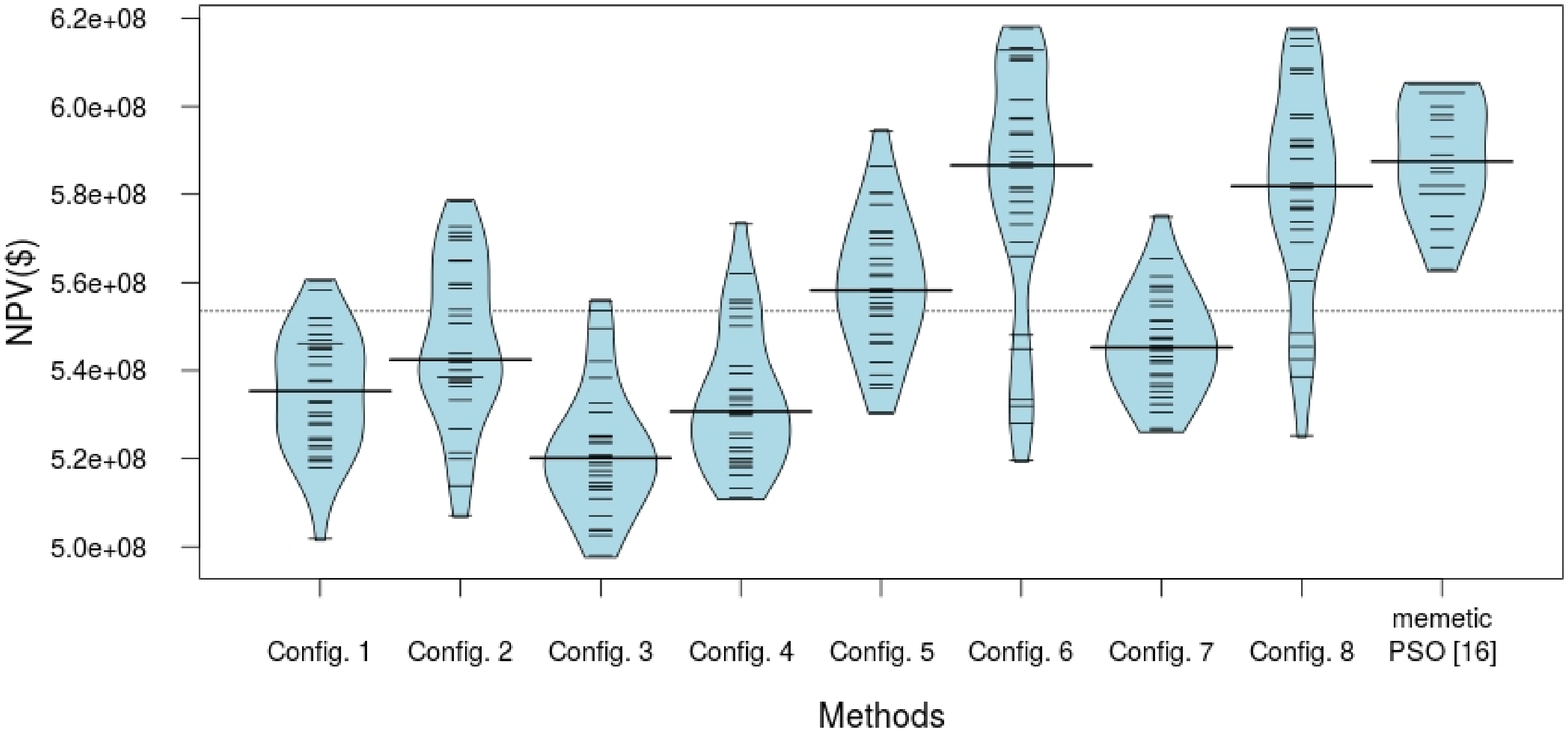} 
\caption{\label{bean3}Bean plots of the NPV distribution for Case~3. The short horizontal lines are the NPVs for all trials. The median is represented by the thick horizontal black line. The average of all observations is presented as the dashed horizontal line.}
\end{figure}

Figures~\ref{bean1}, \ref{bean2}, and \ref{bean3} show bean plots of the NPV distribution~(over $30$ trials) for each DE configuration for Cases~1-3, and for memetic PSO~\cite{Humphries2014}. The short horizontal lines are the NPVs for all trials. The median is represented by the thick horizontal black line. The average of all observations is presented as the dashed horizontal line.

\begin{figure}[h]
\centering
\includegraphics[width=\columnwidth]{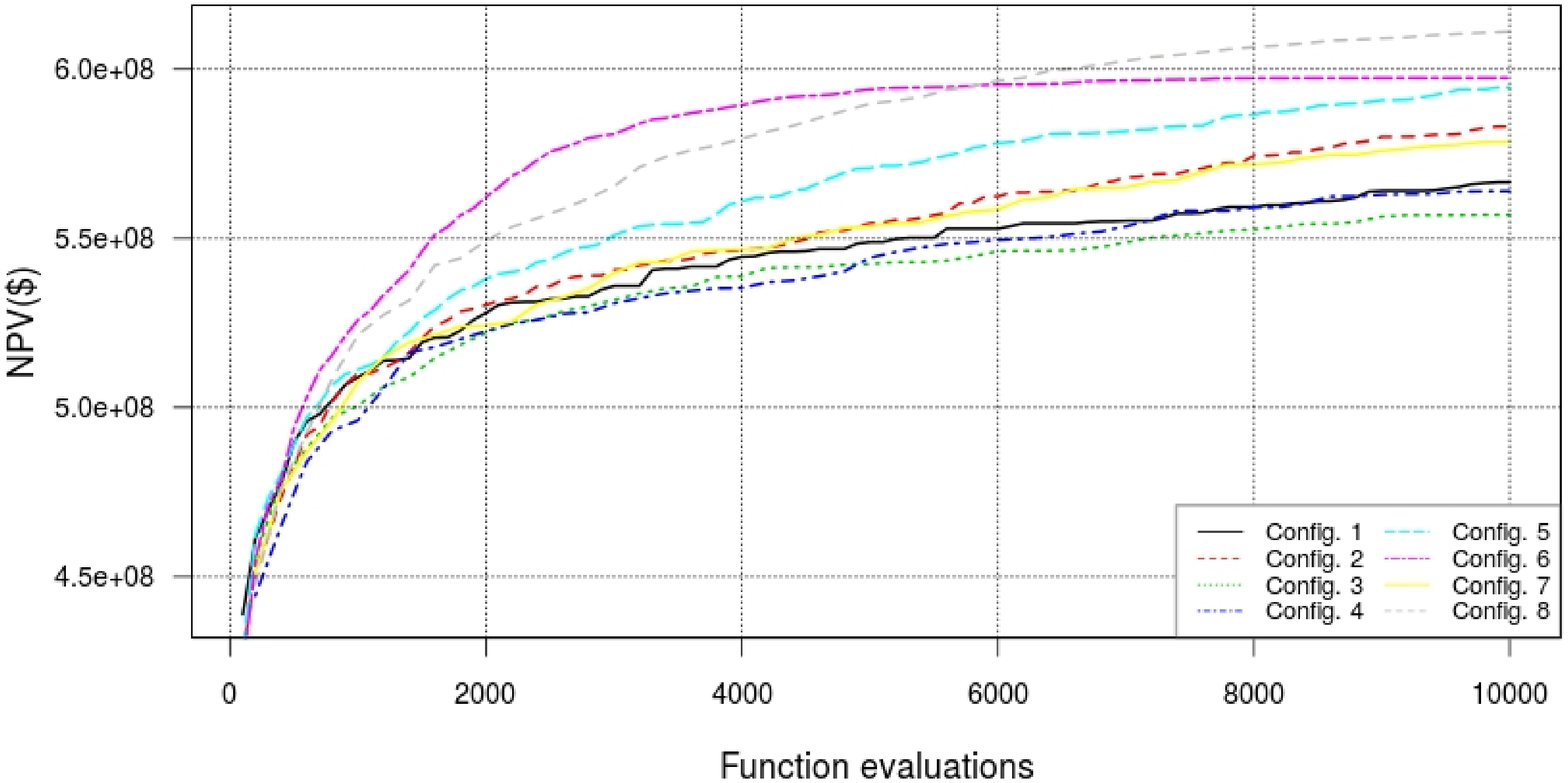} 
\caption{\label{converg1}Convergence plots of the NPV versus the number of function evaluations for Case~1. Each curve represents the mean performance of the $30$ trials for each DE configuration.}
\end{figure}

\begin{figure}[h]
\centering
\includegraphics[width=\columnwidth]{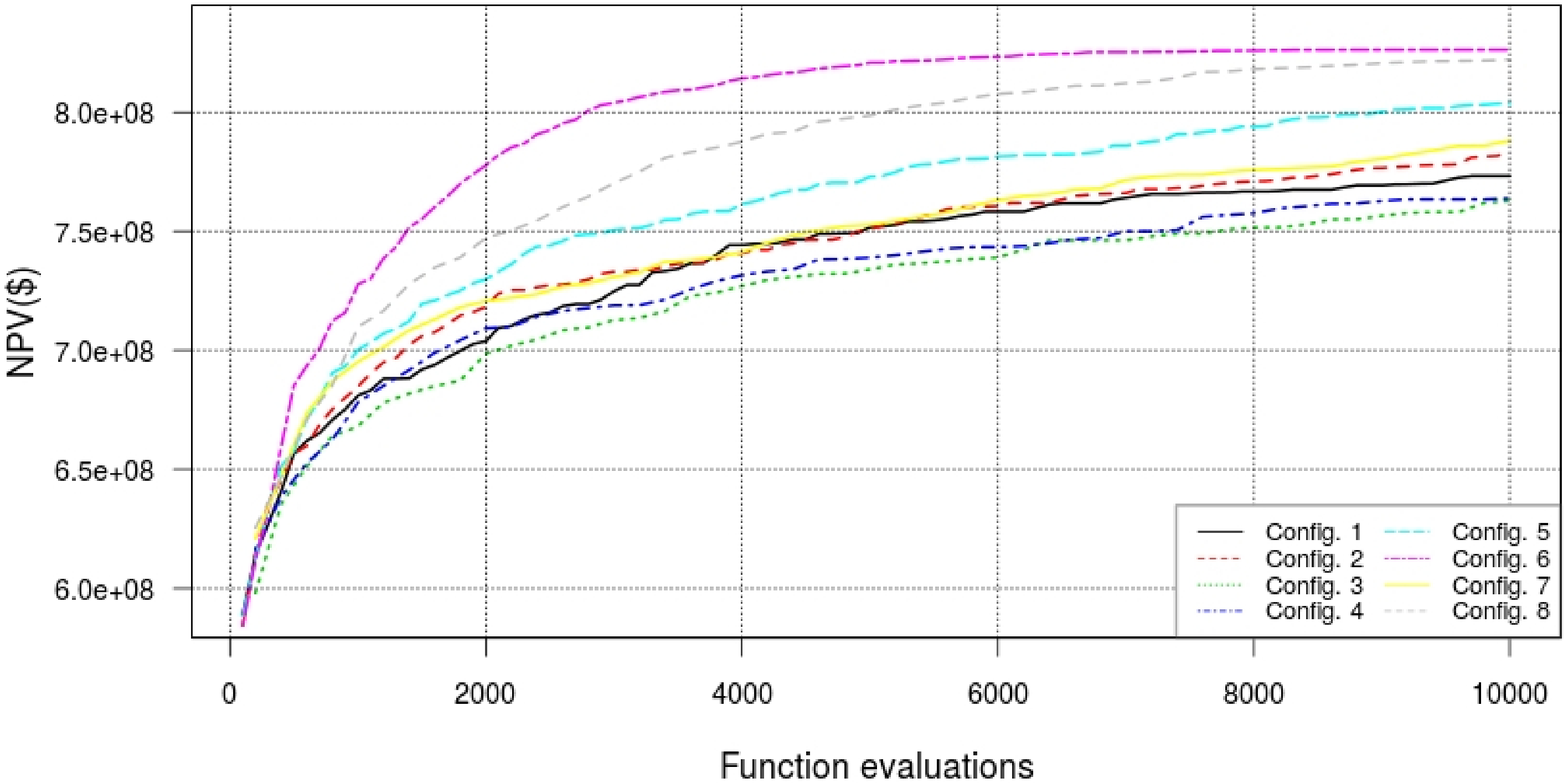} 
\caption{\label{converg2}Convergence plots of the NPV versus the number of function evaluations for Case~2. Each curve represents the mean performance of the $30$ trials for each DE configuration.}
\end{figure}

\begin{figure}[h]
\centering
\includegraphics[width=\columnwidth]{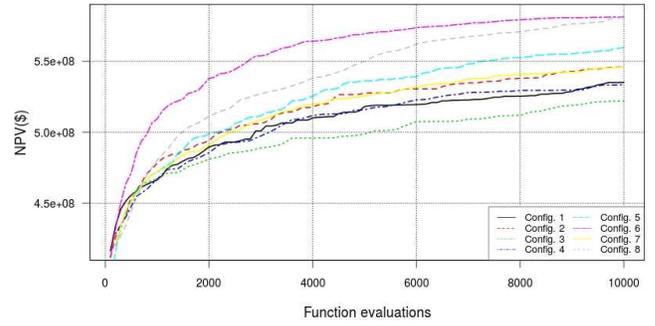} 
\caption{\label{converg3}Convergence plots of the NPV versus the number of function evaluations for Case~3. Each curve represents the mean performance of the $30$ trials for each DE configuration.}
\end{figure}

Figures~\ref{converg1}, \ref{converg2}, and \ref{converg3} show the NPVs with respect to the number of function evaluations for each DE configuration for Cases~1-3. Each curve represents the mean performance of the $30$ trials for each configuration.

\begin{table}[h]
\begin{centering}
\caption{\label{results1}Results for Case~1 described in Section~\ref{case_study}.}
\par\end{centering}
\centering{}%
\scalebox{0.89}{
\begin{tabular}{|c|c|c|c|c|c|}
\cline{2-6} 
\multicolumn{1}{c|}{} & \multicolumn{5}{c|}{\textbf{NPV ($\$ \times 10^8$)}}\tabularnewline
\hline 
\textbf{Method} & \textbf{Best} & \textbf{Worst} & \textbf{Mean} & \textbf{SD} & \textbf{Median}\tabularnewline
\hline 
Config. $1$ & $5.99$ & $5.50$ & $5.67$ & $0.13$ & $5.63$ \tabularnewline
\hline 
Config. $2$ & $6.21$ & $5.56$ & $5.83$ & $0.18$ & $5.82$ \tabularnewline
\hline 
Config. $3$ & $5.82$ & $5.33$ & $5.57$ & $0.12$ & $5.56$ \tabularnewline
\hline 
Config. $4$ & $5.88$ & $5.36$ & $5.64$ & $0.16$ & $5.60$ \tabularnewline
\hline 
Config. $5$ & $6.24$ & $5.61$ & $5.95$ & $0.16$ & $5.96$ \tabularnewline
\hline
Config. $6$ & $6.51$ & $5.34$ & $5.98$ & $0.29$ & $6.03$ \tabularnewline
\hline 
Config. $7$ & $6.06$ & $5.43$ & $5.79$ & $0.15$ & $5.77$ \tabularnewline
\hline 
Config. $8$ & $6.49$ & $5.40$ & $6.11$ & $0.27$ & $6.18$ \tabularnewline
\hline  
memetic PSO \cite{Humphries2014} & $6.30$ & $5.97$ & $6.15$ & - & - \tabularnewline
\hline 
\end{tabular}}
\end{table}

\begin{table}[h]
\begin{centering}
\caption{\label{results2}Results for Case~2 described in Section~\ref{case_study}.}
\par\end{centering}
\centering{}%
\scalebox{0.89}{
\begin{tabular}{|c|c|c|c|c|c|}
\cline{2-6} 
\multicolumn{1}{c|}{} & \multicolumn{5}{c|}{\textbf{NPV ($\$ \times 10^8$)}}\tabularnewline
\hline 
\textbf{Method} & \textbf{Best} & \textbf{Worst} & \textbf{Mean} & \textbf{SD} & \textbf{Median} \tabularnewline
\hline 
Config. $1$ & $8.19$ & $7.39$ & $7.73$ & $0.18$ & $7.69$ \tabularnewline
\hline 
Config. $2$ & $8.29$ & $7.48$ & $7.83$ & $0.19$ & $7.79$ \tabularnewline
\hline 
Config. $3$ & $7.96$ & $7.24$ & $7.63$ & $0.18$ & $7.62$ \tabularnewline
\hline 
Config. $4$ & $8.07$ & $7.40$ & $7.64$ & $0.18$ & $7.59$ \tabularnewline
\hline 
Config. $5$ & $8.40$ & $7.56$ & $8.04$ & $0.20$ & $8.02$ \tabularnewline
\hline
Config. $6$ & $8.64$ & $7.29$ & $8.26$ & $0.29$ & $8.22$ \tabularnewline
\hline 
Config. $7$ & $8.28$ & $7.59$ & $7.88$ & $0.19$ & $7.89$ \tabularnewline
\hline 
Config. $8$ & $8.66$ & $7.14$ & $8.22$ & $0.27$ & $8.18$  \tabularnewline
\hline  
memetic PSO~\cite{Humphries2014} & $8.64$ & $8.04$ & $8.35$ & - & - \tabularnewline
\hline 
\end{tabular}}
\end{table}

\begin{table}[h]
\begin{centering}
\caption{\label{results3}Results for Case~3 described in Section~\ref{case_study}.}
\par\end{centering}
\centering{}%
\scalebox{0.89}{
\begin{tabular}{|c|c|c|c|c|c|}
\cline{2-6} 
\multicolumn{1}{c|}{} & \multicolumn{5}{c|}{\textbf{NPV ($\$ \times 10^8$)}}\tabularnewline
\hline 
\textbf{Method} & \textbf{Best} & \textbf{Worst} & \textbf{Mean} & \textbf{SD} & \textbf{Median} \tabularnewline
\hline 
Config. $1$ & $5.60$ & $5.02$ & $5.35$ & $0.14$ & $5.35$ \tabularnewline
\hline 
Config. $2$ & $5.78$ & $5.07$ & $5.46$ & $0.18$ & $5.42$ \tabularnewline
\hline 
Config. $3$ & $5.56$ & $4.98$ & $5.22$ & $0.15$ & $5.20$ \tabularnewline
\hline 
Config. $4$ & $5.73$ & $5.11$ & $5.34$ & $0.16$ & $5.31$ \tabularnewline
\hline 
Config. $5$ & $5.94$ & $5.31$ & $5.59$ & $0.16$ & $5.58$ \tabularnewline
\hline
Config. $6$ & $6.18$ & $5.19$ & $5.81$ & $0.28$ & $5.87$ \tabularnewline
\hline 
Config. $7$ & $5.75$ & $5.26$ & $5.46$ & $0.12$ & $5.45$ \tabularnewline
\hline 
Config. $8$ & $6.17$ & $5.25$ & $5.81$ & $0.24$ & $5.82$ \tabularnewline
\hline  
memetic PSO~\cite{Humphries2014} & $6.05$ & $5.63$ & $5.89$ & - & - \tabularnewline
\hline 
\end{tabular}}
\end{table}

~\\
\subsection{Discussions}

Considering the experiments in this paper for Case~1, the best NPV, according to Table~\ref{results1}, was found by Config.~6, while the worst one was obtained by Config.~3. Configuration~8 presented the highest median NPV. In turn, Config.~6 showed the highest SD. DE was able to find better solutions than those found in~\cite{Humphries2014} by the decoupled method, but the mean NPV was a little lower, and the worst solutions were lower than those from~\cite{Humphries2014} because of the outliers. Given that we did not employ any local search method, the results of DE are very satisfactory. Moreover, one usually selects the best solution to solve a problem, not the mean solution.

As one can notice in Fig.~\ref{bean1}, the distributions of the final NPVs for each DE configuration are clear\-ly not Gaussian; thus the mean and SD are not adequate measures to describe the distributions. The average of all observations is presented as the dashed horizontal line right above $5.8\times 10^8$.

Configs.~2, 4, 6, and 8, that used $CR=0.9$ and $F=0.5$, presented more variance in the distribution than the other configurations. However, the best results were obtained by Config.~6 and Config.~8. We noticed that Config.~8 gave the highest median, represented by the thick horizontal black line. Configuration~6 obtained the second best result, followed by Config.~5. The worst performance was presented by Config.~3. This information is also shown in Table~\ref{results1}, but the bean plot allows for the visualization of the concentration of solutions with similar values and the identification of outliers. For instance, Config.~6 and Config.~8 found the best solutions, but also some very poor ones. Configuration~8 had a high concentration of solutions close to $6.3\times 10^8$ and $6.0\times 10^8$. On the other hand, based on the plot one may say that the best solution found by Config.~6 was an outlier, since the high concentration was around $6.0\times 10^8$.

Another pattern identified in Fig.~\ref{bean1} is that for \textit{rand/1}, the bigger $N$, the worse the result (see Config.~1-4). For the \textit{current-to-best/1} strategy, the increase in $N$ considerably improved the solutions for Config.~8, but worsened the solutions for Config.~7. We believe that Config.~7 was worse than Config.~5 because when $CR=0.5$, only $50\%$ of the mutated vector will be incorporated into the trial solution. As $N$ for Config.~7 is twice as in Config.~5, it is performing half the number of iterations. Therefore, it would need much more function evaluations to transfer useful information to the offspring.

Regarding the results from~\cite{Humphries2014} shown in Fig.~\ref{bean1}, one can notice that the distributions of the final NPVs seems Gaussian, and presents a high concentration of solutions around $6.17\times 10^8$. In turn, Config.~8 obtained the highest median NPV and a greater concentration of solutions above the median compared to that from~\cite{Humphries2014}, even presenting more variance in the distribution. As shown in Table~\ref{results1} and Fig.~\ref{bean1}, Configs.~6 and 8 gave higher NPVs than the highest NPV found in~\cite{Humphries2014}. Memetic PSO was more stable than DE for this case, giving smaller variance in the results. On the other hand, DE found several solutions considerably better than those of memetic PSO. Therefore, for Case~1 we conclude that even though the memetic PSO presents a higher mean NPV than DE, DE has a higher probability of finding better solutions than the memetic PSO.

In the convergence plots of Fig.~\ref{converg1}, one may observe that Config.~6 was the fastest method to find high quality solutions; it presented the highest NPV using $6,000$ function evaluations, but then stagnated. After $6,000$ function evaluations, Config.~8 obtained the highest NPV. A possible explanation is that these configurations have $90\%$ chance of considering the parent's information ($CR=0.9$). Configuration~6 lost diversity faster than Config.~8, which kept improving because of a bigger $N$. Config.~5 shows a NPV close to Config.~6 in the end, and would probably outperform it if more iterations were allowed. The worst mean performance was presented by Config.~3. The other configurations showed an intermediate performance.

Regarding Case~2, one can see in Table~\ref{results2} that the best solution was found by Config.~8, whereas both Config.~6 and memetic PSO~\cite{Humphries2014} found the second-best solution. On the other hand, the mean value of the simulations in~\cite{Humphries2014} is higher, as is the worst solution found. This means that DE got stuck in local optima more often than the memetic PSO in~\cite{Humphries2014}, or that DE could not properly refine the solutions. Nevertheless, it is reasonable to think that a DE with local search could be more competitive.

When evaluating the bean plots shown in Fig.~\ref{bean2}, we can see that the worst values for Configs.~6 and~8 were outliers. Thus, the worst, mean, and SD from Table~\ref{results2} must be carefully evaluated; they are descriptive statistics of the data, but do not tell the whole story.

Patterns observed in Case~1 bean plots for DE are also valid in Case~2. For instance, \textit{current-to-best/1} was noticeably the best mutation strategy, even though it gave a higher variance. Choosing $CR$ to be larger than $F$ results in better median solutions, except for Config.~4. Also, for \textit{rand/1} a larger $N$ showed worse results than a smaller $N$. The main difference here is that Config.~6 was better than Config.~8, with more solutions near the best ones, i.e, with a high concentration of solutions close to $8.6\times 10^8$, possibly due to a smaller $N$ leading to faster convergence.

Another interesting fact is that the bean plot shapes of Config.~6 and~\cite{Humphries2014} are similar, except for the tail of the bean plot for Config.~6 which is caused by the $4$ poor solutions found. As in Case~1, memetic PSO presented a smaller variance in the results than DE. On the other hand, DE Config.~6 gave more solutions near the best overall even though it obtained a median NPV lower than that from~\cite{Humphries2014}. 

In the convergence plots for Case~2 in Fig.~\ref{converg2}, one can notice that, once again, Config.~6 gave a faster approximation to the best solution, followed by Config.~8. High-quality solutions were found by Config.~6 requiring less than $6,000$ function evaluations, and the remaining function evaluations gave small adjustments in the solution. However, as the DE used in this work does not have a local-search mechanism, the refinement was insufficient to give dramatically better results in most runs. 

For Case~3, one can observe in Table~\ref{results3} that the best NPV was found by Config.~6, and Config.~8 found the second-best NPV. Thus, two DE configurations were able to find better solutions than those found in~\cite{Humphries2014}. However, as in Case~1, the mean performance of DE was a little lower and the worse solutions were lower than obtained in~\cite{Humphries2014}. This confirms that a memetic DE could be more competitive, since DE got stuck in local optima more often than the memetic PSO used in~\cite{Humphries2014}.

Considering the bean plots of DE shown in Fig.~\ref{bean3} for Case~3, Config.~6 presented the highest median, followed by Config.~8. Analogous to Case~1, the worst performance was presented by Config.~3. This information can be seen in Table~\ref{results3} and Fig.~\ref{converg3} as well. As in Cases~1 and~2, the patterns observed in Figs.~\ref{bean1} and~\ref{bean2} are also valid for Fig.~\ref{bean3} of Case~3. For example, worse NPVs are obtained with a large $N$. Better median NPVs are achieved when $CR$ is larger than $F$. Even though presenting a higher variance, \textit{current-to-best/1} was clearly the best mutation strategy. As in Case~2, Config.~6 presented more solutions near to the best ones, i.e, close to $6.18\times 10^8$.

Analogous to Fig.~\ref{bean1}, we can also see in Fig.~\ref{bean3} that Config.~6 and Config.~8 presented higher NPVs than the best NPV found by~\cite{Humphries2014}. Although Config.~6 gave a higher variance in the results it obtained a median NPV close to that from~\cite{Humphries2014}. 

According to Fig.~\ref{converg3}, once again, Config.~6 presented the faster convergence, followed by Config.~8. As mentioned in Section~\ref{DE}, for the configurations that use \textit{current-to-best/1} strategy, smaller $N$ values give faster convergence. Also, Config.~3 showed the worst mean performance. As mentioned for Cases~1 and~2, we believe that all solutions could be better if a local-search was employed since the results from basic DE configurations were competitive with those from memetic PSO~\cite{Humphries2014}.
~\\
\section{Conclusions and Future work}

In this paper, we presented a study of the performance of eight DE configurations in solving a well placement problem. DE was tested with different values for $N$, $CR$, and $F$, and also two well-known distinct mutation strategies. The SPE10 base case was the benchmark dataset, and the simulator used was MRST. We solved three optimization problems with or without constraints, and with or without discounting rate. For each one, a bare control procedure was assumed with a control interval of $2$ years, and the parameters representing spatial coordinates of vertical wells, injectors and producers, were optimized. Performance of DE in solving that problem was compared with results of the memetic PSO from \cite{Humphries2014}. 

According to the preliminary analysis performed in this paper, \textit{current-to-best/1} was visibly the best mutation strategy, in spite of a higher variance in the NPV results. Choosing $CR=0.9$ resulted in better median solutions than with $CR=0.5$. Moreover, a larger $N$ showed worse results. Thus, for the cases tested in this paper one should not expect good results using \textit{rand/1} as mutation strategy or using $CR=0.5$. 

For all three optimization problems considered, DE found better NPVs than those found in~\cite{Humphries2014}. On the other hand, the memetic PSO used in~\cite{Humphries2014} was more stable with the highest mean NPV, and the poorest results found by memetic PSO were higher than the poorest results found by DE.

As future work we intend to use a DE algorithm capable of automatically adapting its parameters, and also add a local-search mechanism. Such a DE variant may not only be as stable as the memetic PSO, but also reach better median values for all problems.

\section*{Acknowledgments}

The authors would like to thank the anonymous reviewers for their valuable comments and suggestions to improve the quality of the paper, and National Council for Scientific and Technological Development (CNPq), grant 248571/2013-3, and NSERC by support this research.
~\\
\bibliographystyle{abbrv}
\bibliography{bibliography/proposal}

\begin{thebibliography}{10}

\bibitem{Afshari2014}
S.~Afshari, M.~R. Pishvaie, and B.~Aminshahidy.
\newblock Well placement optimization using a particle swarm optimization
  algorithm, a novel approach.
\newblock {\em PETROL SCI TECHNOL}, 32(2):170--179, 2014.

\bibitem{Awotunde2014a}
A.~A. Awotunde and N.~Sibaweihi.
\newblock Consideration of voidage-replacement ratio in well-placement
  optimization.
\newblock {\em SPE Economics \& Management}, 6(1), 2014.

\bibitem{Bangerth2006}
W.~Bangerth, H.~Klie, M.~Wheeler, P.~Stoffa, and M.~Sen.
\newblock On optimization algorithms for the reservoir oil well placement
  problem.
\newblock {\em COMPUTAT GEOSCI}, 10(3):303--319, 2006.

\bibitem{Boussaid2013}
I.~Boussa\"{i}d, J.~Lepagnot, and P.~Siarry.
\newblock A survey on optimization metaheuristics.
\newblock {\em INFORM SCIENCES}, 237(0):82 -- 117, 2013.

\bibitem{Bouzarkouna2012}
Z.~Bouzarkouna, D.~Ding, and A.~Auger.
\newblock Well placement optimization with the covariance matrix adaptation
  evolution strategy and meta-models.
\newblock {\em COMPUTAT GEOSCI}, 16(1):75--92, 2012.

\bibitem{Brest2010}
J.~Brest, A.~Zamuda, I.~Fister, and M.~Mau\v{c}ec.
\newblock Large scale global optimization using self-adaptive differential
  evolution algorithm.
\newblock In {\em IEEE C EVOL COMPUTAT}, pages 1--8, July 2010.

\bibitem{Chakraborty2008}
U.~K. Chakraborty.
\newblock {\em Advances in Differential Evolution}.
\newblock Springer, 2008.

\bibitem{Christie2001}
M.~A. Christie and M.~J. Blunt.
\newblock Tenth spe comparative solution project: A comparison of upscaling
  techniques.
\newblock In {\em SPE Reservoir Simulation Symposium}. SPE, Feb. 2001.

\bibitem{Clark2003}
R.~A. Clark and B.~Ludolph.
\newblock Voidage replacement ratio calculations in retrograde condensate to
  volatile oil reservoirs undergoing eor processes.
\newblock In {\em SPE Annual Technical Conference and Exhibition}. SPE, Oct.
  2003.

\bibitem{Das2011}
S.~Das and P.~N. Suganthan.
\newblock Differential evolution: A survey of the state-of-the-art.
\newblock {\em IEEE T EVOLUT COMPUT}, 15(1):4--31, 2011.

\bibitem{Melo2013a}
V.~V. de~Melo and G.~L.~C. Carosio.
\newblock Investigating multi-view differential evolution for solving
  constrained engineering design problems.
\newblock {\em EXPERT SYST APPL}, 40(9):3370--3377, July 2013.

\bibitem{Ding2008}
Y.~Ding.
\newblock Optimization of well placement using evolutionary methods.
\newblock In {\em Europec/EAGE Conference and Exhibition}. SPE, June 2008.

\bibitem{duPlessis2012}
M.~du~Plessis and A.~Engelbrecht.
\newblock Using competitive population evaluation in a differential evolution
  algorithm for dynamic environments.
\newblock {\em EUR J OPER RES}, 218(1):7--20, 2012.

\bibitem{Fonseca2013}
R.~Fonseca, O.~Leeuwenburgh, P.~Van~den Hof, and J.~Jansen.
\newblock Improving the ensemble optimization method through covariance matrix
  adaptation (cma-enopt).
\newblock In {\em Reservoir Simulation Symposium}. SPE, 2013.

\bibitem{Ahmed2014}
A.~Hashim~Ahmed, A.~Awotunde, O.~Mutrif~Siddig, and M.~Jamal.
\newblock A pareto-based well placement optimization.
\newblock In {\em 76th EAGE Conference and Exhibition 2014}. EAGE, June 2014.

\bibitem{Humphries2014}
T.~D. Humphries, R.~D. Haynes, and L.~A. James.
\newblock Simultaneous and sequential approaches to joint optimization of well
  placement and control.
\newblock {\em COMPUTAT GEOSCI}, 18(3-4):433--448, 2014.

\bibitem{Isebor2013}
O.~J. Isebor, L.~J. Durlofsky, and D.~Echeverr\'{i}a~Ciaurri.
\newblock A derivative-free methodology with local and global search for the
  constrained joint optimization of well locations and controls.
\newblock {\em COMPUTAT GEOSCI}, pages 1--20, 2013.

\bibitem{Islam2012}
S.~Islam, S.~Das, S.~Ghosh, S.~Roy, and P.~Suganthan.
\newblock An adaptive differential evolution algorithm with novel mutation and
  crossover strategies for global numerical optimization.
\newblock {\em IEEE T SYST MAN CY B}, 42(2):482--500, Apr. 2012.

\bibitem{Leskinen2009}
J.~Leskinen, F.~Neri, and P.~Neittaanm\"{a}ki.
\newblock Memetic variation local search vs. life-time learning in electrical
  impedance tomography.
\newblock In M.~Giacobini, A.~Brabazon, S.~Cagnoni, G.~Di~Caro, A.~Ek\'{a}rt,
  A.~Esparcia-Alc\'{a}zar, M.~Farooq, A.~Fink, and P.~Machado, editors, {\em
  Applications of Evolutionary Computing}, volume 5484 of {\em LECT NOTES
  COMPUT SC}, pages 615--624. Springer, 2009.

\bibitem{Lie2012}
K.~Lie, S.~Krogstad, I.~S. Ligaarden, J.~R. Natvig, H.~M. Nilsen, and
  B.~Skaflestad.
\newblock Open-source matlab implementation of consistent discretisations on
  complex grids.
\newblock {\em COMPUTAT GEOSCI}, 16(2):297--322, 2012.

\bibitem{Lyons2013}
J.~Lyons and H.~Nasrabadi.
\newblock Well placement optimization under time-dependent uncertainty using an
  ensemble kalman filter and a genetic algorithm.
\newblock {\em J PETROL SCI ENG}, 109(0):70 -- 79, 2013.

\bibitem{Mendes2005}
R.~Mendes and A.~S. Mohais.
\newblock Dynde: a differential evolution for dynamic optimization problems.
\newblock In {\em Congress on Evolutionary Computation}, pages 2808--2815,
  2005.

\bibitem{Mezura-Montes2008}
E.~Mezura-Montes, M.~Reyes-Sierra, and C.~Coello.
\newblock Multi-objective optimization using differential evolution: A survey
  of the state-of-the-art.
\newblock In U.~Chakraborty, editor, {\em Advances in Differential Evolution},
  volume 143 of {\em STUD COMP INTELL}, pages 173--196. Springer, 2008.

\bibitem{Morales2010}
A.~N. Morales, H.~Nasrabadi, and D.~Zhu.
\newblock A modified genetic algorithm for horizontal well placement
  optimization in gas condensate reservoirs.
\newblock In {\em SPE Annual Technical Conference and Exhibition}. SPE, Sept.
  2010.

\bibitem{Neri2011}
F.~Neri, G.~Iacca, and E.~Mininno.
\newblock Disturbed exploitation compact differential evolution for limited
  memory optimization problems.
\newblock {\em INFORM SCIENCES}, 181(12):2469 -- 2487, 2011.

\bibitem{Neri2010}
F.~Neri and V.~Tirronen.
\newblock Recent advances in differential evolution: a survey and experimental
  analysis.
\newblock {\em ARTIF INTELL REV}, 33(1-2):61--106, 2010.

\bibitem{Nwankwor2013}
E.~Nwankwor, A.~Nagar, and D.~Reid.
\newblock Hybrid differential evolution and particle swarm optimization for
  optimal well placement.
\newblock {\em COMPUTAT GEOSCI}, 17(2):249--268, 2013.

\bibitem{Onwunalu2010}
J.~E. Onwunalu and L.~J. Durlofsky.
\newblock Application of a particle swarm optimization algorithm for
  determining optimum well location and type.
\newblock {\em COMPUTAT GEOSCI}, 14(1):183--198, 2010.

\bibitem{Peaceman1978}
D.~Peaceman.
\newblock Interpretation of well-block pressures in numerical reservoir
  simulation.
\newblock {\em SPE Journal}, 18(3):183--194, 1978.

\bibitem{Sarma2008}
P.~Sarma and W.~H. Chen.
\newblock Efficient well placement optimization with gradient-based algorithms
  and adjoint models.
\newblock In {\em Reservoir Simulation Symposium}. SPE, Feb. 2008.

\bibitem{Storn1997}
R.~Storn and K.~Price.
\newblock {Differential Evolution - A Simple and Efficient Heuristic for Global
  Optimization over Continuous Spaces}.
\newblock {\em J GLOBAL OPTIM}, 11(4):341--359, Dec. 1997.

\bibitem{Talbi2009}
E.-G. Talbi.
\newblock {\em Metaheuristics - From Design to Implementation}.
\newblock Wiley, 2009.

\bibitem{Temizel2014}
C.~Temizel, S.~Purwar, K.~Urrutia, A.~Abdullayev, F.~Md~Adnan, A.~Agarwal,
  A.~Garcia, and S.~E. Gorucu.
\newblock Optimization of well placement in real-time production optimization
  of intelligent fields with use of local and global methods.
\newblock In {\em SPE Annual Technical Conference and Exhibition}. SPE, Oct.
  2014.

\bibitem{Tirronen2008}
V.~Tirronen, F.~Neri, T.~K\"{a}rkk\"{a}inen, K.~Majava, and T.~Rossi.
\newblock An enhanced memetic differential evolution in filter design for
  defect detection in paper production.
\newblock {\em EVOL COMPUT}, 16(4):529--555, Dec. 2008.

\bibitem{Yeten2003}
B.~Yeten, L.~J. Durlofsky, and K.~Aziz.
\newblock Optimization of nonconventional well type, location, and trajectory.
\newblock {\em SPE Journal}, 8(3):200--210, 2003.

\bibitem{Zhang2009}
J.~Zhang and A.~C. Sanderson.
\newblock Jade: Adaptive differential evolution with optional external archive.
\newblock {\em IEEE T EVOLUT COMPUT}, 13(5):945--958, 2009.

\bibitem{Zhang2010}
K.~Zhang, G.~Li, A.~C. Reynolds, J.~Yao, and L.~Zhang.
\newblock Optimal well placement using an adjoint gradient.
\newblock {\em J PETROL SCI ENG}, 73(3-4):220--226, 2010.

\bibitem{Zhou2012}
X.~Zhou, Z.~Wu, and H.~Wang.
\newblock Elite opposition-based differential evolution for solving large-scale
  optimization problems and its implementation on gpu.
\newblock In {\em Parallel and Distributed Computing, Applications and
  Technologies (PDCAT), 2012 13th International Conference on}, pages 727--732,
  Dec. 2012.

\bibitem{Zou2013}
D.~Zou, J.~Wu, L.~Gao, and S.~Li.
\newblock A modified differential evolution algorithm for unconstrained
  optimization problems.
\newblock {\em Neurocomputing}, 120(0):469 -- 481, 2013.
\newblock Image Feature Detection and Description.

\end{thebibliography}

\end{document}